\documentstyle[12pt]{article}     

\title{\mbox{}      }
\author{\mbox {}}
\begin{document}
\begin{bf}
\begin{Large}
\begin{center}
Extremely Non-symmetric, Non-multiplicative, Non-commutative Operator Spaces  \\
\end{center}
\end{Large}
\end{bf}

\begin{large}
\begin{center}
Waclaw Szymanski \\

\end{center}
\end{large}

\textbf {Abstract}. Motivated by importance of operator spaces contained in the set of all scalar multiples of isometries ($MI$-spaces) in a separable Hilbert space for $C^*$-algebras and E-semigroups we exhibit more properties of such spaces. For example, if an $MI$-space contains an isometry with shift part of finite multiplicity, then it is one-dimensional. We propose a simple model of a unilateral shift of arbitrary multiplicity and show that each separable subspace of a Hilbert space is the range of a shift. Also, we show that $MI$-spaces are non-symmetric, very unfriendly to multiplication, and prove a Commutator Identity which elucidates the extreme non-commutativity of these spaces.    

\section*{1. Genesis and Justification}
$B(H)$ is the algebra of all linear, bounded operators in a separable Hilbert space $H$ over the complex numbers $\mathbf{C}$, with the identity $I$.
A \emph {subspace} of $H$ is always closed. A \emph {shift} is always unilateral. An \emph {operator space} is a linear subspace of $B(H)$. For $A_1,...,A_n \in B(H)$ $span(A_1,...,A_n)$ is the set of all linear combinations of $A_1,...,A_n$.

Operator spaces contained in the set $MI$ of all scalar multiples of all isometries in a Hilbert space are the subject of investigation in this paper. For brevity, an operator space contained in $MI$ will be called an $MI-space$. Another possible name: "a subspace of $MI$" is misleading because $MI$ is not a linear space. As will be shown, these are strange spaces, indeed. Just for a start, the shift of multiplicity one and its square cannot belong to one $MI$-space, because their sum does not belong to $MI$. This example will be further explained after Proposition 3.2.. Even though the set $MI$ is a semigroup with operator multiplication, it turns out - cf. Proposition 3.2.- that on $MI$-spaces this multiplicative structure trivializes.

My interest in $MI$-spaces came from an attempt to extend the following result of H. Radjavi and P. Rosenthal [R-R]:
\emph{Each linear space of operators contained in the set of all normal operators is commutative}. With John B. Conway [Con-Sz] we replaced "normal" by "hyponormal" in that theorem and showed that such result is false. Trying to understand what really went "wrong" for hyponormals, in [Cat-Sz] $MI$-spaces were introduced (without that name) and proved to be the culprit. Corollary 3.3. of [Cat-Sz] reads: \emph{If $C$ is a class of operators that contains MI then there are $A,B\in C$ such that $span(A,B)\subset MI \subset C$ and $A,B$ do not commute}. Since the class of hyponormal operators contains $MI$, the hyponormal case followed. Even though the attempt to extend the above [R-R] result failed (so far), $MI$-spaces appeared. Concerning their commutative properties - well - they are really bad - worse than found in [Cat-Sz] - cf Corollary 3.5. This is the first justification why $MI$-spaces are worth attention. I call it the \emph {operator theory justification}.

Another justification comes from $C^*$-algebras. This connection was already made in [Cat-Sz]. A \emph{Cuntz algebra} $O_n$ is a universal $C^*$-algebra generated by isometries $S_1,...S_n \in B(H)$ such that $S_1S_1^*+...+S_nS_n^* = I$ - cf. [D]. In [Cat-Sz] Corollary 2.6. states that if $S_1,...S_n \in B(H)$ are generators of a Cuntz algebra then $span(S_1,...S_n)$ is an $MI$-space. Also a slight generalization of the converse of this result is proved there. I call ths the \emph{$C^*$-algebra justification}.

One more justification comes from continuous tensor product systems of Hilbert spaces introduced by William Arveson in [A1] as a continuous analogue of Fock spaces. It turns out that such product systems are a basic structure in studying semigoups of endomorphisms of $B(H)$ called E-semigroups - cf [A2]. Proposition 2.1. of [A1]
says that if $\alpha$ is a non-zero normal *-endomorphism of $B(H)$ then there are isomerties $V_1,V_2,...$ with mutually orthogonal ranges such that $\alpha(A) = \sum V_n AV_n^*$ for each $A \in B(H)$. The linear space $\mathcal {E}$$ =\{T \in B(H): \alpha(A)T=TA$ for each $A \in B(H)\}$ is norm closed and $T^*S$ commutes with $B(H)$ for each $T,S \in $ $\mathcal {E}$ , therefore    
$T^*S$ is a scalar multiple of $I$. By Proposition 2.1. in the next section, $ \mathcal{E}$ is an $MI$-space. For a concrete product system and a semigroup $\alpha_t$, $t \geq 0$, of normal *-automorphisms of $B(H)$ the operator spaces
$\mathcal{E}$$ _t$ defined as above for $\alpha _t$ play a fundamental role in E-semigroup theory. I call this justification the \emph {E-semigroup justification}. 

In summary, $MI$-spaces appear naturally in three areas: operator theory, $C^*$-algebras, and E-semigroups.

\section*{2. Geometry}
In this section geometric aspects of $MI$-spaces will be discussed. In particular, a geometric model of a shift will be presented in Theorem 2.6., which is, perhaps, of interest on its own.

The basic result on $MI$-spaces is

\vspace {7 mm}

\textbf {2.1. Theorem} ([Cat-Sz, Theorem 2.3]). \emph {Suppose $\mathcal{S}$ $\subset B(H)$ is a linear space. Then $\mathcal{S}$ is an $MI$-space if and only if for each $A,B \in$ $\mathcal{S}$ there is $\lambda \in$ $ \mathbf{C}$ such that $B^*A=\lambda I$.}

\vspace {7 mm}

Consider the mapping $<,>_0 :B(H)\times B(H) \rightarrow B(H)$ defined by $<A,B>_0 =B^*A$ for $A,B \in B(H)$. This mapping satisfies all defining conditions of an inner product , except, in general, being scalar-valued. Denote by $\mathbf {C}$$I$ all scalar multilples of $I$. Theorem 2.1 can be restated as

\vspace{7 mm}

\textbf {2.2. Theorem} ([Cat-Sz, Theorem 2.4]). \emph {Suppose $\mathcal{S}$ $\subset B(H)$ is a linear space. Then $\mathcal{S}$ is an $MI$-space if and only if the restriction of $<,>_0$ to $\mathcal{S}$$\times$$\mathcal{S}$ takes  values in $\mathbf {C}$$I$.}

\vspace {7 mm}

Therefore on an $MI$-space $\mathcal{S}$ we introduce the inner product as follows: given $A,B\in$ $\mathcal{S}$ we let $<A,B>=\lambda$ such that $B^*A=\lambda I$ from Theorem 2.1., that is,  

$<A,B>_0= B^*A = <A,B>I$.

The norm defined by this inner product is the same as the operator norm in $B(H)$ because $\|A\|^2 = \|A^*A\|=<A,A>$
for $A \in$ $\mathcal{S}$.

\vspace {7 mm}

\textbf {2.3. Proposition}. \emph{The norm closure of an MI-space is an MI-space}

\textbf {Proof}. Let $\mathcal{S}$ be an $MI$-space. Take a sequence $A_n \in $ $\mathcal{S}$ and $A \in B(H)$.
Then there are sequences $\lambda _n \in \mathbf {C}$ and $V_n \in B(H)$ isometries such that $A_n = \lambda_n V_n$.
Suppose $A_n \rightarrow A$. Then $A_n^*A_n \rightarrow A^*A $ and  $A_n^*A_n=|\lambda_n|^2 V_n^*V_n = |\lambda_n|^2 I$.  
Therefore $ |\lambda_n|^2$ converges to a non-negative number $|\lambda|^2$ for some $\lambda \in \mathbf{C}$.
Hence $A^*A=|\lambda|^2I$. By [Cat-Sz, Proposition 2.1], $A \in MI$. q.e.d.

\vspace {7 mm}

Therefore, the closure of an $MI$-space is a Hilbert space. Notice that the last proof works for any subset of $MI$.

Let $\mathcal{S}$ be an $MI$-space. Two elements $A,B \in$ $\mathcal{S}$ are orthogonal if $B^*A=0$, which means that $A,B$ have orthogonal ranges. Thus orthonormal vectors in $\mathcal{S}$ are isometries with mutually orthogonal ranges. 

In the $C^*$-algebra justification the isometries $S_1,...S_n$ form an orthonormal basis of $span(S_1,...S_n)$. In the E-semigroup justification , it is proved in [A1, Proposition 2.1] that $V_1, V_2, ...$ are the orthonormal basis of $ \mathcal{E}$. 

\vspace {7mm}

According to the celebrated Wold decomposition theorem, for each isometry $A \in B(H)$ there are unique subspaces $H_u$, $H_s$ of $H$ which reduce $A$ such that $H = H_u \oplus H_s$, the part $A_u$ in $H_u$ is unitary, the part $A_s$ in $H_s$
is a shift, and $A=A_u \oplus A_s$. Therefore, the range of $A$ is $AH=H_u \oplus A_s H_s$, thus, $H \ominus AH = H_s \ominus A_s H_s$ is the wandering space for the shift part $A_s$ of $A$. From this we get immediately the following generalization of Proposition 2.10 in [Cat-Sz]:

\vspace {7mm}

\textbf {2.4. Proposition}. \emph {If an $MI$-space $\mathcal{S}$ contains an isometry $A$ whose shift part has finite multiplicity then $\mathcal{S}$ =$span(A)$, thus dim $\mathcal{S}$ = 1}.

\textbf {Proof}. Take $B$ in the orthogonal complement of $A$ in the Hilbert space $cl$ $\mathcal{S}$ = the norm closure of
$\mathcal{S}$, that is, $A^*B = 0$. This is justified by Proposition 2.3.. Then $BH \subset ker A^* = H \ominus AH =$ the wandering space of $A_s$ - cf. remark above. Since $A_s$ has finite multiplicity, $dim H \ominus AH$ is finite. Since, by Proposition 2.3., $B \in MI$, this is possible only if $B=0$. Thus $\mathcal{S}$ $\subset$ $cl$ $\mathcal{S}$ = $span(A)$.  q.e.d.

\vspace {7 mm}

Therefore, if $\mathcal{S}$ is an $MI$-space and $dim$ $\mathcal{S}$ $> 1$ then the shift part of each isometry in $\mathcal{S}$ has infinite multiplicity. In particular,

\vspace {7mm}

\textbf {2.5. Corollary}. \emph {Suppose $\mathcal{S}$ is an $MI$-space.}

\emph {a. If $\mathcal{S}$ contains a unitary operator $A$ then $\mathcal{S}$ = $span(A)$}.

\emph {b. If $\mathcal{S}$ contains $I$ then $\mathcal{S}$ = $\mathbf {C}$$I$.}

\vspace {7 mm}

The remarks before the last proposition  show also that what really matters when considering orthonormal systems in $MI$-spaces is the shift parts if the isometries involved. Therefore, now we turn to shifts. The next proposition is elementary. It is included here for the sake of completeness. Let $\mathbf{N}$ $=\{1,2,3,...\}$.

\vspace {7 mm}

\textbf {2.6. Proposition}. \emph {Suppose $E = \{e_n: n \in \mathbf{N}\}$ is an orthonormal basis of $H$.}
\emph {If $A_0 :E \rightarrow H$ is a mapping whose range consists of orthonormal vectors then there is a unique isometry $A \in B(H)$ such that $A|E = A_0$}.

\textbf {Proof}. Take any $x \in H$. Since $E$ is an orthonormal basis of $H$, there are unique $\alpha_n \in$ $\mathbf {C}$
such that $x = \sum \alpha_n e_n$ and  $\sum | \alpha _n |^2 < \infty$. Define $Ax=\sum \alpha_n A_0 e_n$. The mapping $A:H \rightarrow H$ is well-defined and preserves inner product. Linearity follows by a standard argument. It is plain that $A|E = A_0$ and that such $A$ is unique.  q.e.d.

\vspace {7 mm}

Now a shift with a given wandering space will be constructed. The construction relies on the following remarkable property
of countable sets: \emph {If $X$ is a finite or countable set then $X \times \mathbf {N}$ is countable}.

\vspace {7 mm}

\textbf {2.7. Theorem}. \emph {For each subspace $M$ of $H$ with infinite dimensional $H \ominus M$ there is a shift for which $M$ is the wandering space.}

\textbf {Proof}. Suppose $dim M = m$ is finite or countable. Choose an orthonormal basis $e_{10}, e_{20}, ...,e_{m0}$ of $M$.
Let $X=\{1,...,m\}$. Choose an orthonormal basis of $H \ominus M$ indexed by $X \times \mathbf {N}$ as follows:

\[
\begin{array}{ccc}
e_{11}, e_{21},...,e_{m1}\\
e_{12}, e_{22},...,e_{m2}\\
........................\end{array}
\]

This is possible because $X \times \mathbf {N}$ is countable. Then $E=\{e_{jk}:(j,k)\in X \times (\mathbf {N}\cup\{0\})\}$ is an orthonormal basis of $H$. Define $A_0(e_{jk})=e_{j,k+1}$ for $(j,k) \in X \times (\mathbf {N}\cup \{0\})$. By Proposition 2.6., there is a unique isometry $A\in B(H)$ such that $Ae_{jk}=A_0e_{jk}=e_{j,k+1}$ for $(j,k) \in X \times (\mathbf {N}\cup \{0\})$. Since $A^p M=span(e_{1p},...,e_{mp})$ for $p \in \mathbf {N}\cup \{0\}$, the subspaces $A^p M$
and $A^q M$ are mutually orthogonal for $p,q \in \mathbf {N}\cup \{0\}$, $p \not=q$. Moreover, $H$ is the orthogonal sum of all $A^p M$, $p \in \mathbf {N}\cup \{0\}$, because $E$ is an orthonormal basis of $H$. Hence $A$ is a shift with wandering space $M$.  q.e.d.

\vspace {7 mm}

The construction of the shift in the above proof provides us with a very simple, yet useful model.

\vspace {7 mm}

\textbf {2.8. Corollary}. \emph {Each infinite dimensional subspace of a Hilbert space is the range of a shift.}

\textbf {Proof}. Suppose $K$ is an infinite dimensional subspace of $H$. By Theorem 2.7., there is a shift $A$ with wandering space $H \ominus K$. Since $H \ominus AH$ is the wandering space for $A$, we conclude $K=AH$.  q.e.d. 

\vspace {7 mm}

In the operator theory and $C^*$-algebra justification $MI$-spaces are finite dimensional.
In the E-semigroup justification to avoid trivial cases the $MI$-spaces $\mathcal E$$_t$ have to be infinite dimensional. 

Corollary 2.8. shows, in particular, how to construct $MI$-spaces with any dimension and prescribed ranges of isometries in their orthonormal bases. To get an $MI$-space $\mathcal{S}$ with $dim$$\mathcal{S}$ $=d$ finite or countable and mutually orthogonal ranges $K_1,...,K_d$ of isometries in the orthonormal basis of $\mathcal{S}$ just use Corollary 2.8. to get shifts $A_1,...,A_d$ with desired properties.

Finally, suppose an $MI$-space should have an orthonormal basis consisting of isometries, some of which with non-trivial unitary part. This can be done exaclty the same way as described above for shifts, using the following

\vspace {7 mm}

\textbf {2.9. Corollary}. \emph{Suppose $K$ is a subspace of $H$, $K_u$ is a subspace of $K$ with infinite dimensional $K \ominus K_u$, and $U \in B(K_u)$ is a unitary operator. Then there is an isometry $A \in B(H)$ with unitary part $U$ and range $K$}.

\textbf {Proof}. Let $A_s$ be the shift in $H \ominus K_u$ with range $K \ominus K_u$ as constructed in Theorem 2.7.
and Corollary 2.8.. Let $A=U\oplus A_s$ on $H = K_u \oplus (H \ominus K_u)$. Since $AH = UK_u \oplus A_s(H \ominus K_u) = K_u \oplus (K \ominus K_u) = K$, $A$ satisfies all requirements.  q.e.d.

\vspace {7mm}

Therefore, not only the range, but also the unitary part of an isometry can be arbitrarily prescribed. The only restriction is $K_u \subset K$, but Wold decomposition makes it necessary.

\section* {3. Algebra}

Throughout this section $\mathcal S$ is an $MI$-space. Now we will justify properties of $MI$-spaces in the title of this paper. First, symmetry. 

\vspace {7 mm}

\textbf {3.1. Proposition}. 

\emph {a) If $A\in MI$ is such that $A^* \in MI$ then $A$ is a scalar multiple of a unitary operator.}

\emph {b) If $A \in$ $\mathcal S$ is such that $A \not= 0$ and $A^* \in$ $\mathcal S$ then $\mathcal S$ = $span(A)$.}

\textbf {Proof}.a) Suppose $A \not= 0$. Since $A$ and $A^*$ are in $MI$, there are $\lambda, \mu \in \mathbf {C}$
both non-zero, and isometries $V,W \in B(H)$ such that $A = \lambda V$, $A^* = \mu W$. Hence $\overline{\lambda}V^*= \mu W$ and $(\overline{\lambda}/ \mu)V^* = W$ is an isometry. Therefore, $ker V^* = 0$. Thus $VH=H$ and $V$ is unitary. 
 
b) If such $A$ exists then , by a), it is a non-zero scalar multiple of a unitary operator. By Corollary 2.5.a), $\mathcal S$ = $span(A)$.   q.e.d.

\vspace {7 mm}

Therefore, if an $MI$-space is more than one dimensional then it cannot contain the adjoint of any of its elements. Now, let us turn to multiplicative properties.

\vspace {7 mm}

\textbf {3.2. Proposition}. 

\emph{a) If $A \in$ $\mathcal S$ then $\{B \in B(H): AB \in$ $\mathcal S$$\}$ = $\mathbf{C}$$I$.}

\emph{b) Suppose $\mathcal S$ $\not=$ $\mathbf{C}$$I$ and $A,B \in$ $\mathcal S$. Then $AB \in$ $\mathcal S$
         if and only if $B=0$}.

\emph {c) If $\mathcal S$ $\not=$ $\mathbf{C}$$I$ and $A, A^k \in$} $\mathcal S$ \emph{for some $k \in$ $\mathbf N$, $k>1$, then $A=0$}.  

\textbf {Proof}. Suppose $A \in$ $\mathcal S$ and $B \in B(H)$ is such that $AB \in$ $\mathcal S$. By Proposition 2.1.,
$A^*(AB) = \lambda I$ for some $\lambda \in \mathbf{C}$. But $A^*A = \mu I$ for some $\mu \in \mathbf{C}$. Therefore $B \in$ $\mathbf{C}$$I$. This proves a). Now, if $A,B \in$ $\mathcal S$ are such that $AB \in$ $\mathcal S$ then, by a), $B= \lambda I$ for some $\lambda \in$ $\mathbf{C}$. If $\mathcal S$ $\not=$ $\mathbf{C}$$I$ then, by Corollary 2.5.b), $I \not\in$ $\mathcal S$. Therefore, $\lambda = 0$ and $B = 0$, which proves b). Part c) follows from b).   q.e.d.

\vspace {7 mm}

Part c) of this proposition explains the example with the shift of multiplicity one and its square given at the beginning of Section 1. Now it is plain that $MI$-spaces are very unfriendly to the operator multiplication. A simple way of thinking about an operator space could be considering an operator algebra and forgetting about multplication. Not here.
As we see, for an $MI$-space there is no non-trivial chance even to contain a power of its element, not to mention the algebra of polynomials in its element. 

Finally, we come back to where we started in the operator theory justification, but with a broader perspective.

If $H$ is just any Hilbert space with inner product $(,)$ then the number $\|x\|^2 \|y\|^2 - |(x,y)|^2$ for $x,y \in H$ does not seem to have any particular significance. It is, certainly, non-negative due to Schwarz inequality. In $MI$-spaces, however, this number seems to be rather important - it turns out to be the "measure of non-commutativity" for operators in such spaces, as the following proposition shows. For operators $A,B \in B(H)$ 
\emph{the commutator} is defined by $[A,B]=AB-BA$.

\vspace {7 mm}

\textbf {3.3. Theorem}. \emph {If $A,B \in$ $\mathcal S$ then $[A,B] \in MI$ and}

\begin {center}
$[A,B]^*[A,B] = 2 (\|A\|^2 \|B\|^2 -  |<A,B>|^2)I$
\end {center}

\textbf {Proof}. Just compute:

$[A,B]^*[A,B]= (AB-BA)^*(AB-BA)= B^*A^*AB-B^*A^*BA-A^*B^*AB+A^*B^*BA=$

$B^*\|A\|^2B - B^*<B,A>A - A^*<A,B>B + A^*\|B\|^2 A =$ 

$(\|A\|^2 \|B\|^2 - <B,A><A,B> - <A,B><B,A> + \|B\|^2\|A\|^2)I=$ 

$(2\|A\|^2\|B\|^2 - 2 |<A,B>|^2)I$.

Proposition 2.1. of [Cat-Sz] implies that $[A,B] \in MI$.  q.e.d.

\vspace {7 mm}

It seems the formula proved in this theorem deserves a name - I propose to call it the \emph {Commutator Identity}.
How non-commutative are orthonormal vectors in $\mathcal S$ ? The Commutator Identity gives the answer :

\vspace {7 mm}

\textbf {3.4. Corollary}. \emph {If $A,B \in$ $\mathcal S$ are orthonormal then $[A,B]^*[A,B]= 2I$ and $\|[A,B]\| = \sqrt 2$}.

\vspace{7 mm}

Can we say anything "positive" about commuting in $MI$-spaces ? The answer is a definite no. $MI$-spaces are the most noncommutative linear spaces encountered in the world of operator theory.

\vspace {7 mm}

\textbf {3.5. Corollary}. \emph { Suppose $\mathcal S$ is an $MI$-space.}

\emph {a) $A,B \in$ $\mathcal S$ commute if and only if they are linearly dependent.}

\emph {b) $\mathcal S$ is commutative if and only if $dim$$\mathcal S$ = 0 or $dim$$\mathcal S$ = 1.}

\textbf {Proof}. That $A,B$ commute means $[A,B] = 0$, which is the same as $[A,B]^*[A,B]=0$. By the Commutator Identity, this is equivalent to equality in the Schwarz inequality for $A,B$. As every inner product child knows, this is equivalent to linear dependence of $A,B$.
This proves a). If $\mathcal S$ is commutative and $dim$$\mathcal S$ $\not= 0$ choose a non-zero $A \in$ $\mathcal S$.  Then each $B \in$ $\mathcal S$ commutes with $A$. By a), $A,B$ are linearly dependent. Thus $dim$$\mathcal S$ = 1. q.e.d.

\vspace {10 mm}
\begin{center}
{\bf REFERENCES}
\end{center}

[A1] W. Arveson - \emph {Continuous Analogues of Fock Spaces} Memoirs Amer. Math. Soc. vol 80, no 409 - 1989.

[A2] W. Arveson - \emph {Noncommutative Dynamics and E - semigroups} - Springer Verlag - 2003.

[Cat-Sz] X. Catepillan, W. Szymanski - \emph {Linear Combinations of Isometries} - Rocky Mountains J. Math., 34 (2004), 187-193.

[Con-Sz] J. B. Conway, W. Szymanski - \emph{Linear Combinations of Hyponormal Operators} - Rocky Mountains J. Math., 18 (1988), 695-705.

[D] K. Davidson -\emph {$C^*$-algebras by Example} - AMS Fields Institute Monographs no 6, 1996.  

[R-R] H. Radjavi, P. Rosenthal - \emph {On Invariant Subspaces and Reflexive Algebras} - Amer. J. Math. 91 (1969), 638-692.

\vspace {20 mm}

Waclaw Szymanski - Department of Mathematics, West Chester University,West Chester PA 19382

wszymanski@wcupa.edu

\end{document}